\documentclass[12pt,oneside,english]{amsart}
\usepackage[latin1]{inputenc}
\usepackage{amssymb}

\makeatletter
\newtheorem{theorem}{Theorem}

\newtheorem{proposition}{Proposition}
\newtheorem{example}{Example}
\newtheorem{definition}{Definition}
\newtheorem{remark}{Remark}

\newtheorem{corollary}{Corollary}

\usepackage{babel}

\makeatother
\begin{document}
\baselineskip=17pt

\title{Overpseudoprimes, Mersenne numbers and Wieferich primes}

\author{Vladimir Shevelev}
\address{Department of Mathematics \\Ben-Gurion University of the
 Negev\\Beer-Sheva 84105, Israel. e-mail:shevelev@bgu.ac.il}


\begin{abstract}
 We introduce a new class of pseudoprimes - so-called ``overpseudoprimes" which is a special subclass of super-Poulet
 pseudoprimes. Denoting via $h(n)$ the multiplicative order of 2 modulo $n$, we show that odd number $n$ is overpseudoprime if and only if the value of h(n) is invariant of all divisors $d>1$ of $n$. In particular, we prove that all composite Mersenne numbers $2^{p}-1$, where $p$ is prime, and squares of Wieferich primes are overpseudoprimes.

\end{abstract}

\maketitle

\section{Introduction}

     Sometimes the numbers $M_n=2^n-1, \enskip n=1,2,\ldots $, are called Mersenne numbers, although this name is
usually reserved for numbers of the form

\begin{equation}\label{1}
M_p= 2^p -1
\end{equation}

where $p$ is prime. In our paper we use the latter name. In this form numbers $M_p$ at the first time were studied by Marin Mersenne (1588-1648) at least in 1644 (see in [1, p.9] and a large bibliography there).\newline

           We start with the following  simple observation. Let $n$ be odd and $h(n)$ denote the multiplicative order of 2 modulo $n$.

\begin{theorem}\label{t1} Odd $d>1$ is a divisor of $M_p$ if and only if $h(d) =p.$
\end{theorem}
\bfseries Proof. \mdseries If $d>1$ is a divisor of $2^{p}-1$, then $h(d)$ divides prime $p$. But $h(d)>1$.
Thus, $h(d)=p$. The converse statement is evident.$\blacksquare$
\begin{remark}\label{r1}This observation for prime divisors of $M_p$ belongs to Max Alekseyev $($ see his
comment to sequence $A122094$ in \cite{5}$).$
\end{remark}
\newpage
     In our paper, we in a natural way  introduce a new class  $\mathbb{S}$  of  pseudoprimes  and show that it contains those and only those odd numbers $n$ for which $h(n)$ is invariant of all divisors $d>1$ of $n$. In particular, it contains all composite Mersenne numbers
     and, at least, squares of all Wieferich primes \cite{6}. We also give a generalization to arbitrary base $a>1.$

\section{A class of pseudoprimes}

For an odd $n >1$, consider the number $r=r(n)$ of distinct cyclotomic cosets of 2 modulo $n$  [2, pp.104-105]. E.g., $r(15)=4$
since for $n=15$ we have the following 4 cyclotomic cosets of 2: $\{1,2,4,8\}, \{3,6,12,9\},\{5,10\},\newline\{7,14,13,11\}$.

Note that, if $C_1,\ldots, C_r$ are all different cyclotomic cosets of 2$\mod n$, then

\begin{equation}\label{2}
\bigcup^r_{j=1}C_j=\{1,2,\ldots, n-1\},\qquad C_{j_1}\cap C_{j_2}=\varnothing , \;\; j_1\neq j_2.
\end{equation}

For the least common  multiple of $|C_1|, \ldots, |C_r|$ we have

\begin{equation}\label{3}
[|C_1|,\ldots,|C_r|]= h(n).
\end {equation}

(This follows easily, e.g., from Exercise 3, p. 104 in \cite{3}).

It is easy to see that for odd prime $n=p$ we have

\begin{equation}\label{4}
|C_1|=\ldots=|C_r|=h(n),
\end {equation}

such that

\begin{equation}\label{5}
p= rh + 1.
\end {equation}

\begin{definition}
We call odd composite number $n$ overpseudoprime to base 2 $(n\in\mathbb{S}_2)$ if

\begin{equation}\label{6}
n=r(n) h(n)+1.
\end {equation}

\end{definition}

Note that

$$
2^{n-1}= 2^{r(n)h(n)}\equiv 1\pmod{n}.
$$

Thus, $\mathbb{S}_2$ is a subclass of Poulet class of pseudoprimes to base $2$ (see[6]).
\newpage
\begin{proposition}
If $n\in\mathbb{S}_2,$ then $(\ref{4})$ satisfies.
\end{proposition}
\bfseries Proof. \mdseries If not all $|C_i|$ are the same, then $n-1=\sum |C_i|<h(n)\sum 1=h(n)r(n),$ which contradicts to the definition of overpseudoprime. $\blacksquare$

\begin{theorem}\label{t2}

Let $n$ be odd composite number with the prime factorization

\begin{equation}\label{7}
n=p_1^{l_1}\cdots p_k^{l_k}.
\end {equation}

Then $n$ is overpseudoprime to base 2 if and only if for all nonzero vectors $(i_1, \ldots, i_k)\leq (l_1, \ldots, l_k)$ we have

\begin{equation}\label{8}
h(n)=h(p_1^{i_1}\cdots p_k^{i_k}).
\end {equation}

\end{theorem}

\bfseries Proof 1. \mdseries   Taking into account Proposition \ref{1}, we can suppose that (\ref{4}) satisfies. It is well known that

$$
\sum_{d|n}\varphi(d)=n,
$$

where $\varphi(n)$ is Euler function.
Thus, by (\ref{7})

\begin{equation}\label{9}
\sum_{0\leq i_j \leq l_j, \enskip j=1,\ldots,k}\varphi(p_1^{i_1}\cdots p_k^{i_k})=n.
\end {equation}

     Consider a fixed nonzero vector $(i_1,\ldots, i_k)$ and numbers not exceeding $n$ of the form

\enskip\begin{equation}\label{10}
m=m(i_1,\ldots, i_k) = ap_1^{l_1-i_1}\cdots p_k^{l_k-i_k}, \enskip a\geq1,\enskip (a,p_1^{i_1}\cdots p_k^{i_k} )=1.
\end {equation}

The number of numbers
(\ref{10}) ( or the number of different values of $a\leq \frac{n}{p_1^{l_1-i_1}\cdots p_k^{l_k-i_k}}$) equals to
\begin{equation}\label{11}
\varphi\left(\frac{n}{p_1^{l_1-i_1}\cdots p_k^{l_k-i_k}}\right)=\varphi(p_1^{i_1}\cdots p_k^{i_k}).
\end {equation}

Since for every such $m,$ we have $|C_m|=|C_{ap^{l_1-i_1}...p^{l_k-i_k}}|\leq h(p^{l_1-i_1}...p^{l_k-i_k}),$ then

$$r(m)=\varphi(p_1^{i_1}...p_k^{i_k})/|C_{p^{l_1-i_1}...p^{l_k-i_k}}|\geq$$
\enskip\begin{equation}\label{12}
\varphi(p_1^{i_1}...p_k^{i_k})/h(p^{l_1-i_1}...p^{l_k-i_k}).
\end {equation}

     Thus,
\newpage
$$
r(n)\geq \sum_{0\leq i_j \leq l_j, \enskip j=1,\ldots,k,\enskip not\enskip all\enskip i_j=0} r(m)\geq
$$

\begin{equation}\label{13}
\sum \varphi \left({p_1^{i_1}\cdots p_k^{i_k}} \right)/ h(p_1^{l_1-i_1}\cdots p_k^{l_k-i_k}).
\end {equation}

Note that, by (\ref{4}),

\begin{equation}\label{14}
h(n)\geq  h\left(p_1^{l_1-i_1}\cdots p_k^{l_k-i_k}\right).
\end {equation}

     Thus, by (\ref{13}) and (\ref{9}), we have

\begin{equation}\label{15}
r(n) \geq \frac {1}{h(n)}\sum_{0\leq i_j \leq l_j, \enskip j=1,\ldots,k\enskip not \enskip all \enskip i_j=0}\varphi(p_1^{i_1}\cdots p_k^{i_k})=\frac{n-1}{h(n)},
\end {equation}

and, moreover, the equality holds if and only if for all nonzero vectors $(i_1,\ldots,i_k) \leq (l_1,\ldots, l_k)$, (\ref{8}) is valid. In only this case $r(n)h(n)+1=n$ and $n$ is overpseudoprime. $\blacksquare$\newline\newline

\bfseries Proof 2. \mdseries Suppose that (\ref{4}) holds, such that the cardinality of every coset is $|C(n)|.$ Let $T=T(n)$ be a cyclic group with the generative element $\tau: \tau(i)=2i \pmod n,$ such that $T$ acts on the set $\{1,2,...,n-1\}.$ It is clear that $T$ has order $h(n)$ and the number of its orbits is $r(n)$. Let $d>1$ be a divisor of $n.$ Then $T(d)$ is a cyclic group with the generative element $t: t(i)=2i \pmod d,$ such that $T(d)$ acts on the set $\{1,2,...,d-1\}.$ It is clear that, since $d|n,$ the cardinality of cosets $C(n)$ modulo $n$ equals to the cardinality of cosets $C(d)$ modulo $d$ and the latter equals to $h(d).$ Therefore, $h(d)=h(n).$ $\blacksquare$
\begin{corollary}\label{1} Every two overpseudoprimes to base $2$ $n_1$ and $n_2$ for which $h(n_1)\neq h(n_2)$ are coprimes.
 \end{corollary}

\begin{corollary}\label{2}  Mersenne number $M_p$ is  either prime or overpseudoprime to base $2.$

 \end{corollary}
\bfseries Proof\enskip\mdseries  follows straightforward from Theorems 1-2. $\blacksquare$
\newline     By the definition (see [6]), a Poulet number all of whose divisors $d$ satisfy $d|2^d-2 $ is called a\slshape super-Poulet number.\upshape

\begin{corollary}\label{3}$\mathbb{S}_2$  is a subset of super-Poulet pseuduprimes to base $2.$
\end{corollary}
\bfseries Proof.\enskip\mdseries  Let $n\in\mathbb{S}_2$.\enskip If $1< d|n$ then, by Theorem 2,\enskip $d$\enskip itself is either prime or overpseudiprime to base 2, i.e.
$2^{d-1}\equiv 1\enskip(\mod d). \blacksquare$
\newpage
The following example shows that $\mathbb{S}_2$ is a \slshape proper\upshape\enskip subset of super-Poulet (or strong) pseuduprimes to base $2.$
\begin{example}\label{e1}  Consider a super-Poulet pseudoprime $[5, A001262]$

$$
n=314821 = 13 \cdot 61 \cdot 397.
$$

We have $[5, A002326]$

$$
h(13)=12, \enskip h(61)=60, \enskip h(397) =44.
$$

Thus $n$ is not an overpseudoprime to base $2.$
\end{example}

     Note, that if for primes $p_1 < p_2$ we have $h(p_1)= h(p_2)$  then $h(p_1p_2)= h(p_1)$ and $n=p_1p_2$ is
     overpseudoprime to base 2. Indeed, $h(p_1p_2)\geq h(p_1)$.\enskip But
$$
2^{h(p_1)} = 1+ kp_1 = 1+tp_2.
$$

Thus, $k=sp_2$   and

$$
2^{h(p_1)} = 1+ sp_1p_2.
$$

     Therefore, $h(p_1p_2)\leq h(p_1)$ and we are done. By the same way obtain that if $h(p_1)=\ldots=h(p_k)$ then
$n= p_1 \ldots p_k$ is overpseudoprime to base 2.

\begin{example}\label{e2}  Note that

$$
h(53)= h(157) = h(1613) = 52.
$$
Thus,$$n=53\cdot157\cdot1613=13421773$$
is overpseudoprime to base $2.$
\end{example}

And what is more, by the same way, using Theorem \ref{t2} we obtain the following result.
\begin{theorem}\label{t3}If $ p_i^{l_i}, i=1,\ldots,k$, are overpseudoprimes to base $2$ such that $h(p_1)=\ldots=h(p_k)$ then
$n=p_1^{l_1}\cdots p_k^{l_k}.$ is overpseudoprime to base $2.$
\end{theorem}
\;\;\;\;\;
\section{The $(w+1)$-th power of Wieferich prime of order $w$ is overpseudoprime to base 2}
\begin{definition}\label{d2} A prime\enskip $p$ is called a Wieferich prime $($cf. \cite{6}$)$,\enskip if\enskip
$2^{p-1}\equiv1\pmod {p^2}$; a prime $p$ we call a Wieferich prime of order $w\geq1,$ if\enskip $p^{w+1}\| 2^{p-1}-1$.
\end{definition}
\newpage
\begin{theorem}\label{t4}A prime $p$ is a Wieferich prime of order more or equal to $w,$ if and only if $p^{w+1}$ is overpseudoprime to base $2.$
\end{theorem}

\bfseries Proof.\enskip\mdseries Let prime $p$ be Wieferich prime of order at least $w.$  Let $2^{h(p)}=1+kp$. Note that\enskip$h(p)$\enskip divides $p-1$. Using the condition, we have
$$2^{p-1}-1= (kp+1)^{\frac{p-1}{h(p)}}-1=(kp)^{\frac{p-1}{h(p)}}+\ldots+kp{\frac{p-1}{h(p)}}\equiv0\pmod {p^{w+1}}.$$
Thus, $k\equiv0\pmod {p^{w}}$ and $2^{h(p)}\equiv1\pmod{p^{w+1}}$. Therefore, $h(p^{w+1})\leq h(p)$
and we conclude that
$$ h(p)=h(p^2)=\ldots=h(p^{w+1}).$$
Hence, by Theorem \ref{t2}, $p^{w+1}$ is overpseudoprime to base 2. The converse statement is evident.$\blacksquare$
\begin{theorem}\label{t5}
If overpseudoprime to base $2$ number $n$ is not multiple of square of a Wieferich prime, then $n$ is squarefree.

\end{theorem}

\bfseries Proof. \mdseries   Let $n=p_1^{l_1} \ldots p_k^{l_k}$ and, say, $l_1\geq 2$. If $p_1$ is not a Wieferich prime then $h(p_1^2)$ divides $p_1(p_1-1)$ but does not divide $p_1-1$. Thus, $h(p_1^2)\geq p_1$. Since $h(p_1)\leq p_1-1$ then $h(p_1^2) > h(p_1)$ and by Theorem \ref{2}, $n$ is not overpseudoprime to base 2. $\blacksquare$\newline

The following theorem is a generalization of a known property of Mersenne numbers.

\begin{theorem}\label{t6}
Let $q$ be a prime divisor of $2^{p}-1$ such that $q^2|2^{p}-1$. Then $q^{w}\| 2^{p}-1)$ if and only if $q$ is a Wieferich prime of order $w-1.$

\end{theorem}

\bfseries Proof.\enskip\mdseries  Let $q^{w}|2^{p}-1,\enskip w\geq2.$ Since by Theorem 1,\enskip$ h(q)=p$ then we have
$h(q^{w})\leq h(q)$. Thus, $h(q^{w})=h(q^{w-1})=\ldots=h(q)=p$ and $p$ is a Wieferich prime of order at least $w-1.$
If also $h(q^{w+1})=h(q)=p$ then $2^{p}\equiv1\pmod {q^{w+1}}$ and $q^{w}\nparallel 2^{p}-1)$ .$\blacksquare$

  Note that, an algorithm of search a large prime which as the final result could be not a Mersenne prime is the following: we seek a prime $q$ not exceeding $\sqrt{M_p}$ for which $h(q)=p$; if such prime is  absent, then $M_p$ is prime; if we found a prime $q$, then we seek a prime $q_1\leq \sqrt{\frac{M_p}{q}}$ for which $h(q_1)=p$ and if such prime is absent, then $\frac{M_p}{q}$ is a (large) prime etc.

      Note also that, the problem of the infinity of Mersenne primes is equivalent to the problem of infinity primes $p$ for which the equation $h(x)=p$ has no solutions not exceeding $2^{\frac p 2}$.
At last, notice that, for the only known Wieferich primes  1093 and 3511, we have $h(1093)=364, h(3511)=1755$
\newpage
(see sequence A002326 in \cite{5}). Thus, since 364 and 1755 are not prime, they divide none of Mersenne numbers. The important question is: \slshape  do exist  Wieferich primes $p$ for which $h(p)$ is prime?\upshape \enskip If the conjecture of Guy [1, p.9] about the existence of nonsquarefree Mersenne numbers is true, then we should say ``yes".
\section{ Overpseudoprime to base $a$}
    Here we consider a natural generalization. Let $a$ be integer more than 1. If $(n.a)=1$ denote $h_a(n)$ the multiplicative order of $a$ modulo $n$. Furthermore, denote by $r_a(n)$ the number of cyclotomic cosets of $a$  $\pmod {n}$:  $C_1,\ldots, C_{r_a(n)},$ such that (\ref{2}) satisfies. Let $p$ be a prime which does not divide $a$. It is easy to see that $h_a(p)r_a(p)=p-1.$
  \begin{definition}
We call odd composite number $n$, for which $(n,a)=1$, overpseudoprime to base $a$ $(n\in\mathbb{S}_a)$ if

\begin{equation}\label{16}
n=r_a(n) h_a(n)+1.
\end {equation}

\end{definition}
The following theorem is proved in the same way as Theorem \ref{t2}.
\begin{theorem}\label{t7}

Let $n$ be composite number for which $(n,a)=1$ with the prime factorization

\begin{equation}\label{17}
n=p_1^{l_1}\cdots p_k^{l_k}.
\end {equation}

Then $n$ is overpseudoprime to base $a,$ if and only if for all nonzero vectors $(i_1, \ldots, i_k)\leq (l_1, \ldots, l_k),$ we have

\begin{equation}\label{18}
h_a(n)=h_a(p_1^{i_1}\cdots p_k^{i_k}).
\end {equation}

\end{theorem}
     Furthermore, putting, for a prime $p,$,
      \begin{equation}\label{19}
      M_p^{(a)}=\frac{a^{p}-1}{a-1},
\end {equation}

       we have the following generalization of Theorem \ref{1}.
\begin{theorem}\label{t8} Integer $d>1,$ for which $(d,a(a-1)),$ is a divisor of $M_p^{(a)},$ if and only if $h_a(d) =p.$
\end{theorem}

Thus, from Theorems \ref{t7}, \ref{8} we obtain the following statement.

 \begin{theorem}\label{t9}  If $(M_p^{(a)},a-1)=1,$ then  $M_p^{(a)}$ is  either prime or overpseudoprime of base $a$.
\end{theorem}
\newpage
Thus, at least, if there exist infinitely many composite numbers $ M_p^{(a)}$ with the condition $ (M_p^{(a)}, a-1)=1$, then there exist infinitely many overpseudoprimes of base $a.$
\begin{example} $M_3^{(11)}=133=7\cdot19$ is overpseudoprime of base 11. Indeed, we see that $h_{11}(7)=h_{11}(19)=3$.
\end{example}
\begin{definition}\label{d4} A prime\enskip$p$ is called a \upshape Wieferich prime\slshape \enskip in base $a,$\enskip if\enskip$   a^{p-1}\equiv1\pmod {p^2}$;  a prime\enskip$p$ we call a Wieferich prime in base $a$ of order $w\geq1$ if\enskip$ p^{w+1}\| a^{p-1}-1.$
\end{definition}
\begin{theorem}\label{t10} A prime $p$ is a Wieferich prime in base $a$ of order more or equal to $w,$ if and only if $p^{w+1}$ is overpseudoprime of base $a.$
\end{theorem}

\bfseries Proof\enskip\mdseries is over in the same way as in case of Theorem {t4}.   $\blacksquare$

\begin{example}
$p=5$ is a Wieferich prime in  base $7$ of order $1.$ Thus, 25 is overpseudoprime of base $7.$
\end{example}
Furthemore, we have the following generalization of Theorem {t5}.
\begin{theorem}\label{t11}
If $n$ is overpseudoprime of base $a$ and is not multiple of square of a Wieferich prime, then $n$ is squarefree.

\end{theorem}
       The following result shows that the overpseudoprimes appear not more frequently than the strong pseudoprimes.
\begin{theorem}\label{t12}
If $n$ is overpseudoprime to base $a,$ then $n$ is strong pseudoprime to the same base.
\end{theorem}
\bfseries Proof.\enskip\mdseries  By the definition \cite{6}, if $n$ is odd composite number and $2^{s}\|n-1$ then n is strong pseudoprime to base $a$ in case when either $a^{\frac {n-1} {2^{s}}}\equiv1\pmod n$ or for only $k$, $k=0,\ldots,s-1,$ we have $a^{\frac {n-1} {2^{s-k}}}\equiv-1\pmod n.$
Let $n$ be overpseudoprime to base $a$ such that $2^{t}\|h_a(n)$. Since $h_a(n)|(n-1)$ then $t\leq s$. If
$h_a(n)$ is odd then $t=0$ and $r_a(n)/2^{s}$ is integer. Thus,
$$ a^{\frac {n-1} {2^{s}}}=a^{\frac {h_a(n)r_a(n)} {2^{s}}}\equiv1\pmod n $$
and $n$ is strong pseudoprime to base $a.$
      In case of $t\geq1$ we have
       $$a^{\frac {n-1} {2^{s}}}=a^{\frac {h_a(n)} {2^{t}}\cdot\frac{r_a(n)} {2^{s-t}}}$$
 Put $A=a^{\frac {h_a(n)} {2^{t}}}.$ Note that
$$ (A-1)(A+1)(A^2+1)(A^{2^2}+1)\cdot\ldots\cdot(A^{2^{t-1}}+1)=A^{2^{t}}-1\equiv0\pmod n. $$
Let us show that none divisor $d$ of $n$ divides $A-1.$ Indeed, since $n$ is
\newpage
overpseudoprime to base $a$
then $h_a(d)=h_a(n)$ and the congruence $a^{\frac {h_a(d)} {2^{t}}}\equiv1
\pmod d $\enskip for $t\geq1$
contradicts to the definition of $h_a(d).$ Thus, $(A-1,n)=1$ and we have
$$(A+1)(A^2+1)(A^{2^2}+1)\cdot\ldots\cdot(A^{2^{t-1}}+1)\equiv0\pmod n.$$
Furthermore, none divisor $d$ of $n$ divides the difference $A^{2^{j}}-A^{2^{i}}$ for $ 0\leq i <j \leq t$ because of $(A,n)=1$  and and in view of the impossibility of the congruence
$$ A^{2^{j}-2^{i}}=a^{\frac {h_a(n)} {2^{t}}\cdot(2^{j}-2^{i})}=a^{\frac {h_a(d)} {2^{t}}\cdot(2^{j}-2^{i})}\equiv1\pmod d$$
which, in view of $\frac {2^{j}-2^{i}} {2^{t}}<1,$ contradicts to the definition of $h_a(d).$ Therefore,
$(A^{2^{j}}-A^{2^{i}}, n)=1$ and there exist only $i,i=0.\ldots,t$ such that $A^{2^{i}}\equiv-1\pmod n,$
i.e., $n$ is strong pseudoprime to base $a.$ $ \blacksquare$

\section{ Overpseudoprime to base $a$ as superpseudoprime to the same base}

Quite as in the above, where we proved that every overpseudoprime to base 2 is super-Poulet pseudoprime,
using Theorem 7 it could be proved the following statement.
\begin{theorem}\label{t13}
Every overpseudoprime $n$ to base $a$ is superpseudoprime, i. e., for each divisor $d>1$ of $n,$ we have
\begin{equation}\label{20} a^{d-1}\equiv1\pmod d.
\end{equation}
\end{theorem}
Furthermore, we prove the following.
\begin{theorem}\label{t14} If $n$ is  overpseudoprime to base $a,$ then, for every two divisors $d_1 < d_2$ of  $n$ (including trivial divisors 1 and $n$), we have
\begin{equation}\label{21}
 h_a(n) |\enskip d_2-d_1.
\end{equation}
\end{theorem}
\bfseries Proof.\enskip\mdseries By (\ref{20}), we have $h_a(d_i)$=$h_a(n)|\enskip d_i-1$,\enskip$i=1,2,$ and (\ref{21}) follows. $\blacksquare$
\begin{corollary}\label{c4} If $(M_p^{(a)},a-1)=1,$ then for every two divisors $d_1 < d_2$ of  $M_p^{(a)}$ (including trivial divisors 1 and $M_p^{(a)}$) we have
\begin{equation}\label{22}
 p|\enskip d_2-d_1.
\end{equation}
In particular,this is a property of Mersenne numbers $M_p.$
\end{corollary}
\bfseries Proof\enskip\mdseries follows from Theorems \ref{8} and \ref{14}.$\blacksquare$
\begin{example}\label{e5} For $M_{29}=233\cdot1103\cdot2089,$ we have, in particular, $2089-233=29\cdot64$;\enskip $2089-1103=29\cdot34.$
\end{example}
\newpage
     The following corollary we formulate for Mersenne number $M_r$ although it could be formulated for $M_r^{(a)}$.

\begin{corollary}\label{e5} Let $r$ be prime. Then $M_r$ is prime if and only if the progression $(1+rx)_{x\geq0}$ contains
only prime $p$ with the condition $h_2(p)=r.$
\end{corollary}

\bfseries Proof.\enskip\mdseries If $M_r$ is prime and for some prime $p$ we have $h_2(p)=r$ then, by Theorem 1, $p=M_r.$ Since $1|M_r$ then, by Corollary 4, $r|M_r-1$ i.e. $p=M_r$ is only prime in the progression $(1+rx)_{x\geq0}.$ Conversely, let there exist only prime of the form  $p=1+rx$ with the
condition $h_2(p)=r.$ Then $p|M_r.$ If $M_r$  is composite number, then it is overpseudoprime to base 2 and, for
 another $q|M_r,$ we have $q\equiv p\pmod r$ and, by Theorem 1, $h_2(q)=r.$ This contradicts to the condition. $\blacksquare$

\begin{remark}\label{r2} Till $26.04.08,$ when the author has submitted the sequence
$[5, A 137576]$ under the influence of his paper \cite{4}, he did not touch with the theory of pseudoprimes. He even thought  that the composite numbers $n$ for which $h(n)r(n)=n-1$, probably, do not exist. But after publication of sequence $A137576$ in \cite{5}, Ray Chandler, by the direct calculations, has found a few such numbers. After that the author created a theory which is presented in this paper and found more such numbers of $A141232$ in \cite{5}, using very helpful extended tables of sequences $A002326$ and $A001262$ in \cite{5}, which was composed by T.D.Noe.
\end{remark}

\bfseries Acknowledgment.\mdseries The author is grateful to Max Alekseyev (University of California, San Diego) for useful private correspondence.
\;\;\;\;\;\;\;\;\;\;\;\;\;\;\;


\begin{thebibliography}{5}


\bibitem 1  R.\enskip K. \enskip Guy. \slshape Unsolved Problems in Number Theory, \upshape,2-nd ed. \enskip Springer-Verlag, 1994.
\bibitem 2  F.\enskip J.\enskip  MacWilliams and N.\enskip J.\enskip A.\enskip Sloane, \slshape The Theory of Error-Correcting Codes,\upshape Elsevier/North Holland, 1977.
\bibitem 3  D.\enskip Redmond,\enskip \slshape Number Theory: an Introduction,\enskip\upshape Marcel Dekker,  N.Y., 1996.
\bibitem 4  V.\enskip Shevelev,\enskip\slshape Exact exponent of remainder term of Gelfond's digit theorem in binary case,\upshape \enskip http:// arxiv.org /abs/\enskip 0804.3682.
\bibitem 5   N.\enskip J.\enskip A.\enskip Sloane,\enskip\slshape The On-Line Encyclopedia of Integer Sequences \upshape(http: //www.research.att.com)
\bibitem 6   E.\enskip W.\enskip Weisstein, \enskip\slshape "Poulet number", "Strong pseudoprime", "Wieferich prime", From MathWorld: A Wolfram Web Resource. \upshape(http: //mathworld.wolfram.com/PouletNumber.html)


\end{thebibliography}
\end{document}